\documentclass[colorlinks]{mpi2015-cscpreprint}
\usepackage[american]{babel}
\usepackage{graphicx}
\usepackage[mode=text]{siunitx}
\usepackage{macros}
\begin{document}
\title{Modelling Gas Networks with Compressors: A port-Hamiltonian Approach}
\author[$\ast$]{Thomas Bendokat}%
\author[$\ast$,$\dagger$]{Peter Benner}%
\author[$\ast$]{ \authorcr Sara Grundel}%
\author[$\ast$]{Ashwin S. Nayak}%
\affil[$\ast$]{%
Max Planck Institute for Dynamics of Complex Technical Systems, Magdeburg, Germany.
\email{bendokat@mpi-magdeburg.mpg.de}, \orcid{0000-0002-0671-6291}
\email{benner@mpi-magdeburg.mpg.de}, \orcid{0000-0003-3362-4103}
\email{grundel@mpi-magdeburg.mpg.de}, \orcid{0000-0002-0209-6566}
\email{anayak@mpi-magdeburg.mpg.de}, \orcid{0000-0002-9855-2377}
\vspace{0.3cm}
}
\affil[$\dagger$]{%
Otto von Guericke University, Magdeburg, Germany.
}
%
\shorttitle{Modelling Gas Networks with Compressors}
\shortauthor{T. Bendokat et al.}
\novelty{
  \begin{itemize}
    \item A first attempt at including compressors in port-Hamiltonian formulation for gas networks.
    \item Four nodal models for compressors with a pairwise combination of two specifications and two assumptions.
    \item Numerical implementation for a simple testcase.
    \item Publicly available software material for reproducibility and extension.
  \end{itemize}
}
\abstract{%
  Transient gas network simulations can significantly assist in design and operational aspects of gas networks.
  Models used in these simulations require a detailed framework integrating various models of the network constituents - pipes and compressor stations among others.
  In this context, the port-Hamiltonian modelling framework provides an energy-based modelling approach with a port-based coupling mechanism.
  This study investigates developing compressor models in an integrated isothermal port-Hamiltonian model for gas networks.
  Four different models of compressors are considered and their inclusion in a larger network model is detailed.
  A numerical implementation for a simple testcase is provided to confirm the validity of the proposed model and to highlight their differences.
}
\maketitle
\section{Introduction}

Transportation of hydrogen is a vital strategy in Germany's renewable energy transition.
The federal government has drawn up regulatory policies, convened gas transmission operators and has proposed developing a new core hydrogen-gas network alongside initiatives to repurpose existing pipeline infrastructure.
The initiative brings forth a host of technical challenges related to an optimal design of a gas network and its operation~\cite{MelAP13, Gor2023}.

A gas transportation network comprises mainly of large pipelines interconnected at junctions, but there exist other critical components which can influence gas flow.
Network design and operation decisions can be aided by simulation tools that can predict the transmission operators with vital supply requirement for target demands within the network.
The computational tools serve as an a priori analysis tool into the transient behavior of gas within the network and its components greatly improving efficiency and reliability.
Gas networks crucially contain compressor stations specifically to increase the pressure of incoming gas stream, to compensate for pressure loss due to flow and friction in the pipes.
A significant portion of the operational costs of gas networks are attributed to compressor stations and this differs vastly for hydrogen, natural gas and their blends~\cite{Sema21}.
A digital twin of the entire gas network can aid in such economic forecasts, redundancy planning and emergency response.

An accurate simulation of gas networks requires an integrated model of all its constituents.
Network models have generally emphasized on connected pipeline flow models with only recent studies focusing on including compressor models~\cite{Her07,PamBD16,GyrZ19,BerM21}.
Typical models of compressor stations are based on their power ratings, whereas some models consider it as a node with the pressure scaled by a \textit{compression ratio}~\cite{GyrZ19, BerM21} or as an element scaling the in-flow to a specified \textit{target pressure}~\cite{HimGB21}.
In this study, an integrated model is proposed from an energy-based modelling viewpoint considering that gas networks fit into a port-Hamiltonian framework~\cite{DomHLetal21, HimGB21,HauM21a,HauM21b,HauM23}.
Port-Hamiltonian systems are a combination of port-based modeling and geometric mechanics and are well suited for modeling physical networks such as gas networks.
An introduction and overview of these models can be found in~\cite{VanJ14}.
Such models are suitable to be solved using geometric integration methods, preserving the structure.
An emphasis in this study is placed specifically on including compressor models within the port-Hamiltonian models for gas networks and its software implementation.
\section{Gas Networks}
Gas networks consist of many components: pipes, compressors, resistances, heaters, coolers, valves, etc.~\cite{DomHLetal21}. In this study, we focus on the former two, combining them in a port-Hamiltonian fashion.

\subsection{Pipe Model}

As gas pipelines are long and thin, gas flow in a pipe can be modelled using the \emph{Euler equations} in one spatial dimension~\cite{Her07, DomHLetal21},
\begin{equation}
    \label{eq:Eulerequations}
    \begin{split}
    \partial_t \rho + \partial_x (\rho v) &= 0,\\
    \partial_t(\rho v) + \partial_x(p+\rho v^2) &= - \frac{\lambda}{2D} \rho v \abs{v} - g \rho \sin{\alpha},\\
    \partial_t E + \partial_x ((E+p)v) &= -\frac{k_w}{D} (T-T_w),
    \end{split}
\end{equation}
where \(\rho\) denotes the \emph{density} of the gas, \(v\) denotes the \emph{velocity}, \(p\) is the \emph{pressure}, \(E = \rho(\frac12 v^2 + c_vT +gh)\) is the \emph{total energy} with \emph{specific heat} \(c_v\) and \emph{height} \(h\), \(T\) is the \emph{temperature} of the gas and \(T_w\) is the \emph{temperature of the wall}.
Furthermore, \(\lambda \) is the \emph{pipe friction coefficient}, \(D\) is the \emph{diameter} of the pipe, \(g\) the \emph{gravitational constant}, \(\alpha\) the \emph{inclination angle} of the pipe and \(k_w\) the \emph{heat transfer coefficient} of the wall.
The pressure, density, and temperature are connected via the \emph{equation of state}
\begin{equation}
    \label{eq:stateequation}
    p = \rho R T z(p,T),
\end{equation}
with \emph{gas constant} \(R\) and \emph{compressibility factor} \(z(p,T)\).
The Euler equations~\eqref{eq:Eulerequations} can be transformed into port-Hamiltonian form as they are, but for our case we use the simplified ISO2b model~\cite{DomHLetal21}, by making the following assumptions:
\begin{enumerate}
  \item The model is \emph{isothermal}, meaning the temperature of the pipe is constant in time and space.
  \item The gas velocity \(v\) is insignificant compared to the (constant) speed of sound \(c\), i.e., \(\left({v}/{c}\right)^2\ll 1\).
  \item The pressure in a pipe is related to the density by \(p = c^2 \rho \).
  \item The influence of gravity is neglected.
\end{enumerate}
The ISO2b model then reads
\begin{equation}
    \label{eq:ISO2b}
    \begin{split}
    \partial_t \rho + \partial_x m
        &= 0,\\
    \partial_t m + \partial_x p
        &= - \frac{\lambda}{2D} m \abs{v},
    \end{split}
\end{equation}
with \emph{momentum} \(m = \rho v\).
We apply the ISO2b model in port-Hamiltonian form to a pipe of length \(L\), i.e., \(x \in [0,L]\). To that end, define the \emph{state} \(\z(x,t) = \begin{bmatrix} \rho(x,t)\\ m(x,t) \end{bmatrix}\) and the \emph{effort} \(\e(\z)(x,t) = \begin{bmatrix} p(\rho(x,t))\\ m(x,t) \end{bmatrix}\), as well as
\begin{align*}
    \J(\z) = \begin{bmatrix}
        0 & -D_x \\ -D_x & 0
    \end{bmatrix}\quad \text{and}\quad
    \Rmat(\z) = \begin{bmatrix}
        0 & 0 \\ 0 & \frac{\lambda}{2D} \abs{v}
    \end{bmatrix},
\end{align*}
where \(D_x\) is the differential operator \(D_x f = f_x\) for any function \(f \colon [0,L] \to \R\).
Similar to~\cite{DomHLetal21}, we can enforce the pressure \(p_0\) at the inflow and the momentum \(m_L\) at the outflow of the pipe as boundary conditions, and write the ISO2b model in input-output port-Hamiltonian form as
\begin{equation}
  \label{eq:ISO2b_pH}
  \begin{split}
    \begin{bmatrix}
      I_2 \\
      0
    \end{bmatrix}
      \dot\z
    &=
    \begin{bmatrix}
      \J(\z) - \Rmat(\z) & 0\\
      -U & 0
    \end{bmatrix}
    \begin{bmatrix}
      \e(\z)\\
      Y \e(\z)
    \end{bmatrix}
    +
    \begin{bmatrix}
      0\\
      I_2
    \end{bmatrix}
    \u\\
    \y &= \begin{bmatrix}
      0 & I_2
    \end{bmatrix}
    \begin{bmatrix}
      \e(\z)\\
      Y \e(\z)
    \end{bmatrix},
  \end{split}
\end{equation}
  where
  \begin{align*}
    U \e(\z) = \begin{bmatrix}
      \mathbin{\phantom{-}}\e(\z)_1(0)\\
      {-}\e(\z)_2(L)
    \end{bmatrix},\
    Y\e(\z) = \begin{bmatrix}
      \e(\z)_2(0)\\
      \e(\z)_1(L)
    \end{bmatrix},\
    \u = \begin{bmatrix}
      \mathbin{\phantom{-}}p_0\\
      -m_L
    \end{bmatrix}.
  \end{align*}
Here, \(\e{(\z)}_i(x)\) denotes the \(i\)-th component of \(\e(\z)\), evaluated at \(x\), i.e., for \(x=0\) and \(x=L\), at the inflow and outflow of the pipe, respectively.
\subsection{Compressor Model}
Compressors are used in gas networks to increase gas pressure, which drops along pipes due to friction.
The most commonly used compressors in gas networks are \emph{centrifugal compressors} and \emph{piston compressors}.
While the internal physics of a compressor are quite complex, the relation between the \emph{compression ratio}
\begin{equation}
  \compressionratio = \frac{p_\outlet}{p_\inlet}
\end{equation}
and the \emph{enthalpy} \(H_\ad\) added to the gas is given by~\cite{Her07,WalH17}
\begin{equation}
  H_\ad = z(p_\inlet, T_\inlet) T_\inlet R_s \frac{\kappa}{\kappa-1}\left(\Big(\frac{p_\outlet}{p_\inlet}\Big)^{\frac{\kappa-1}{\kappa}} - 1\right),
\end{equation}
where \(\kappa\) denotes the \emph{isentropic exponent}. For a specific given compressor, the enthalpy \(H_\ad\) is a function of the mass flow rate and the drive speed of the compressor, which is usually approximated from data.
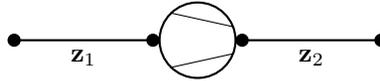
\begin{figure}[ht]
  \begin{center}
	\tikzexternalenable%
	\tikzsetnextfilename{pipe_compressor}%
	\filemodCmp{tikz/pipe_compressor.tex}{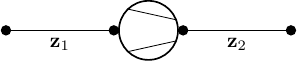}%
	{\tikzset{external/remake next}}{}%
	\def\pipelength{2}

\begin{tikzpicture}
    \draw[thick, {Circle[]}-{Circle[]}]
    (0,0) -- (\pipelength,0) node [pos=0.5, anchor=north] {\(\mathbf{z}_1 \)};
    \node[circle, draw, thick] (circle) at (\pipelength+0.5,0) [minimum size=1cm] {};
    \draw[thin] (circle.135) -- (circle.20);
    \draw[thin] (circle.225) -- (circle.340);
    \draw[thick, {Circle[]}-{Circle[]}] (\pipelength+1,0) -- ++(\pipelength,0) node [pos=0.5, anchor=north] {\(\mathbf{z}_2 \)};
\end{tikzpicture}%
	\tikzexternaldisable%

  \caption{A compressor is inserted between two pipes.}\label{fig:compressor}
  \end{center}
\end{figure}
To insert a compressor between two pipes in a port-Hamiltonian framework, we combine two pipes from~\eqref{eq:ISO2b_pH}, visualized in Fig.~\ref{fig:compressor}, and implement the compressor in the matrix \(\hat\G(\hat\z ) \) and the boundary vector \(\hat\u \) as
\begin{equation}
  \label{eq:compressor}
  \begin{split}
  \hat\E
  \frac{\partial \hat\z}{\partial t}
  &=
  \left(\hat\J(\hat\z )- \hat\Rmat(\hat\z )\right)\hat\e(\hat\z )
  +
  \hat\G\hat\u \\
  \hat\y&= \hat\G^T \hat\e(\hat\z ),
\end{split}
\end{equation}
where
\begin{align*}
  \hat\z &= \begin{bmatrix}
    \z_1\\
    \z_2
  \end{bmatrix},& \hat\e(\hat\z ) &= \begin{bmatrix}
    \e(\z_1)\\
    Y\e(\z_1)\\
    \e(\z_2)\\
    Y\e(\z_2)
  \end{bmatrix},&&\\
  \hat\E &= \begin{bmatrix}
    I_2 & 0  \\
    0 & 0 \\
    0 & I_2\\
    0 & 0
  \end{bmatrix},&
  \hat\J(\hat\z) &= \begin{bmatrix}
    \J_1(\z_1) & 0 & 0 & 0\\
    -U & 0 & 0 & 0\\
    0 & 0 & \J_2(\z_2) & 0\\
    0 & 0 & -U & 0
  \end{bmatrix},&
  \hat\Rmat(\hat\z )&=
  \begin{bmatrix}
   \Rmat_1(\z_1) & 0 & 0 & 0\\
    0 & 0 & 0 & 0\\
    0 & 0 & \Rmat_2(\z_2) & 0\\
    0 & 0 & 0 & 0
  \end{bmatrix}.
\end{align*}
We model compressors as jump conditions under two different assumptions and in two different frameworks, respectively.
\begin{enumerate}
  \item \textbf{Assumption~AV}: Constant velocity of the gas is assumed across the compressor, as in~\cite[Table 9.2]{Els16}.
  In this case, the relationship between the momentum \(m_\inlet \) at the inlet of the compressor and the momentum \(m_\outlet \) at the outlet is given by
  \begin{equation}
    m_\outlet = \tilde c^{\frac{1}{\kappa}} \frac{z(p_\inlet, T_\inlet)}{z(p_\outlet, T_\outlet)} m_\inlet
    \tag{\textbf{AV}}.
  \end{equation}
  In practice, we consider \(z\) to be constant, so that we model \(m_\outlet = \tilde c^{\frac{1}{\kappa}} m_\inlet\).
  \item \textbf{Assumption~AM}: Constant momentum of the gas is assumed across the compressor, as in~\cite{Her07,ZloCB15,SunZ19}. In this case
  \begin{equation}
    m_\outlet = m_\inlet.
    \tag{\textbf{AM}}
  \end{equation}
\end{enumerate}
Under both of those assumptions, we model compressors in the following two frameworks.
\subsubsection{Framework FC}
We choose the compression ratio \(\compressionratio \) of the compressor, meaning the output pressure varies with the input pressure
\begin{equation}
    p_\outlet(\compressionratio,p_\inlet) = \compressionratio p_\inlet.
    \tag{\textbf{FC}}
\end{equation}
In this case, the matrix \(\hat\G(\hat\z )\) in~\eqref{eq:compressor} is given by
\begin{align*}
  \hat\G_{\textbf{FC}}(\hat\z ) = \begin{bmatrix}
    0 & 0 & 0 & 0 \\
    1 & 0 & 0 & 0 \\
    0 & -m_2(0) & 0 & 0 \\
    0 & 0 & 0 & 0\\
    0 & 0 & p_1(L) & 0\\
    0 & 0 & 0 & 1
  \end{bmatrix},
\end{align*}
and, depending on the assumption, the boundary vector \(\hat\u \) is given by
\begin{align*}
  \hat\u^{\textbf{AV}}_{\textbf{FC}} = \begin{bmatrix}
    p_0\\
    \tilde{c}^{-\frac{1}{\kappa}}\\
    \tilde{c}\\
    -m_L
  \end{bmatrix}\quad \text{ or }\quad
  \hat\u^{\textbf{AM}}_{\textbf{FC}} = \begin{bmatrix}
    p_0\\
    1\\
    \tilde{c}\\
    -m_L
  \end{bmatrix}.
\end{align*}
The output \(\hat\y \) in this setting is given by
\begin{align*}
  \hat\y_{\textbf{FC}}(\hat\z ) = \hat\G_{\textbf{FC}}(\hat\z )^T \hat\e(\hat\z ) = \begin{bmatrix}
    m_1(0)\\
    -m_2(0)p_1(L)\\
    p_1(L)m_2(0)\\
    p_2(L)
  \end{bmatrix},
\end{align*}
which leads to external energy exchange
\begin{align*}
  \hat\y_{\textbf{FC}}(\hat\z )^T\hat\u^{\textbf{AV}}_{\textbf{FC}} &= p_1(0)m_1(0) - p_2(L)m_2(L) + \left(\tilde c - \tilde c^{-\frac{1}{\kappa}}\right)p_1(L)m_2(0), \text{ or }\\
  \hat\y_{\textbf{FC}}(\hat\z )^T\hat\u^{\textbf{AM}}_{\textbf{FC}} &= p_1(0)m_1(0) - p_2(L)m_2(L) + \left(\tilde c - 1\right)p_1(L)m_2(0).
\end{align*}
In both cases, we see that for compression ratio \(\tilde c = 1\), the external energy exchange is the energy inserted at the inlet minus the energy removed at the outlet, while for a compression ratio \(\tilde c > 1\), energy is added by the compressor.
\subsubsection{Framework FP}
We choose the output pressure \(p_\outlet \) of the compressor, meaning the compression ratio varies depending on the input pressure
\begin{equation}
  \compressionratio(p_\outlet,p_\inlet) = \frac{p_\outlet}{p_\inlet}.
  \tag{\textbf{FP}}
\end{equation}
In this case, depending on the assumption, the matrix \(\hat\G(\hat\z )\) in~\eqref{eq:compressor} is given by
\begin{align*}
  \hat\G^{\textbf{AV}}_{\textbf{FP}}(\hat\z ) = \begin{bmatrix}
    0 & 0 & 0 & 0 \\
    1 & 0 & 0 & 0 \\
    0 & -m_2(0)p_1(L)^{\frac{1}{\kappa}} & 0 & 0 \\
    0 & 0 & 0 & 0\\
    0 & 0 & 1 & 0\\
    0 & 0 & 0 & 1
  \end{bmatrix}\quad \text{ or }\quad
  \hat\G^{\textbf{AM}}_{\textbf{FP}}(\hat\z ) = \begin{bmatrix}
    0 & 0 & 0 & 0 \\
    1 & 0 & 0 & 0 \\
    0 & -m_2(0) & 0 & 0 \\
    0 & 0 & 0 & 0\\
    0 & 0 & 1 & 0\\
    0 & 0 & 0 & 1
  \end{bmatrix}
\end{align*}
and, depending on the assumption, the boundary vector \(\hat\u \) is given by
\begin{align*}
  \hat\u^{\textbf{AV}}_{\textbf{FP}} = \begin{bmatrix}
    p_0\\
    p_\outlet^{-\frac{1}{\kappa}}\\
    p_\outlet\\
    -m_L
  \end{bmatrix}\quad \text{ or }\quad
  \hat\u^{\textbf{AM}}_{\textbf{FP}} = \begin{bmatrix}
    p_0\\
    1\\
    p_\outlet\\
    -m_L
  \end{bmatrix}.
\end{align*}
The output \(\hat\y \) in this setting, depending on the assumption, is given by
\begin{align*}
  \hat\y^{\textbf{AV}}_{\textbf{FP}}(\hat\z ) = \hat\G^{\textbf{AV}}_{\textbf{FP}}(\hat\z )^T \hat\e(\hat\z ) = \begin{bmatrix}
    m_1(0)\\
    -m_2(0)p_1(L)^{\frac{\kappa+1}{\kappa}}\\
    m_2(0)\\
    p_2(L)
  \end{bmatrix}\quad \text{ or }\quad
  \hat\y^{\textbf{AM}}_{\textbf{FP}}(\hat\z ) = \hat\G^{\textbf{AM}}_{\textbf{FP}}(\hat\z )^T \hat\e(\hat\z ) = \begin{bmatrix}
    m_1(0)\\
    -m_2(0)p_1(L)\\
    m_2(0)\\
    p_2(L)
  \end{bmatrix}
\end{align*}
which leads to external energy exchange
\begin{align*}
  \hat\y^{\textbf{AV}}_{\textbf{FP}}(\hat\z )^T\hat\u^{\textbf{AV}}_{\textbf{FP}} &= p_1(0)m_1(0) - p_2(L)m_2(L) + \left(p_\outlet- p_1(L) \left(\frac{p_1(L)}{p_\outlet}\right)^{\frac{1}{\kappa}}\right) m_2(0), \text{ or }\\
  \hat\y^{\textbf{AM}}_{\textbf{FP}}(\hat\z )^T\hat\u^{\textbf{AM}}_{\textbf{FP}} &= p_1(0)m_1(0) - p_2(L)m_2(L) + \left(p_\outlet - p_1(L)\right)m_2(0).
\end{align*}
In both cases, it is readily seen that if the pressure at the inlet of the compressor \(p_1(L)\) is equal to the output pressure \(p_\outlet \), no energy is added by the compressor.
\subsection{Network Interconnection}
\begin{figure}[h]
  \begin{center}
	\tikzexternalenable%
	\tikzsetnextfilename{pipe_compressor_network}%
	\filemodCmp{tikz/pipe_compressor_network.tex}{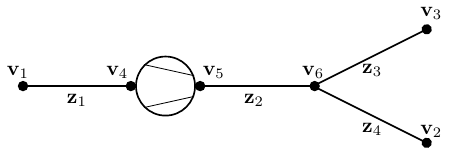}%
	{\tikzset{external/remake next}}{}%
	\def\pipelength{2}

\begin{tikzpicture}
    \draw[thick, {Circle[]}-{Circle[]}]
    (0,0) -- (\pipelength,0) node [pos=0, anchor=south] {\(\mathbf{v}_1 \)} node [pos=0.5, anchor=north] {\(\mathbf{z}_1 \)} node [pos=1, anchor=south east] {\(\mathbf{v}_4 \)};
    \node[circle, draw, thick] (circle) at (\pipelength+0.5,0) [minimum size=1cm] {};
    \draw[thin] (circle.135) -- (circle.20);
    \draw[thin] (circle.225) -- (circle.340);
    \draw[thick, {Circle[]}-{Circle[sep=-3.2pt]}] (\pipelength+1,0) -- ++(\pipelength,0) node [pos=0, anchor=south west] {\(\mathbf{v}_5 \)} node [pos=0.5, anchor=north] {\(\mathbf{z}_2 \)} node [pos=1, anchor=south] {\(\mathbf{v}_6 \)};
    \draw[thick, -{Circle[]}] (\pipelength+\pipelength+1,0) -- ++(\pipelength,1) node [pos=0.5, anchor=north] {\(\mathbf{z}_3 \)} node [pos=1, anchor=south] {\(\mathbf{v}_3 \)};
    \draw[thick, -{Circle[]}] (\pipelength+\pipelength+1,0) -- ++(\pipelength,-1) node [pos=0.5, anchor=north] {\(\mathbf{z}_4 \)} node [pos=1, anchor=south] {\(\mathbf{v}_2 \)};
\end{tikzpicture}%
	\tikzexternaldisable%

  \caption{A simple pipe network with an inserted compressor.}\label{fig:compressor_network}
  \end{center}
\end{figure}
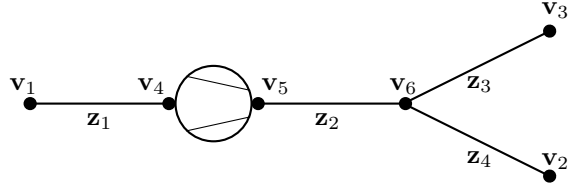

Structure preserving interconnection of the compressor model into a port-Hamiltonian gas network can be done by extension of the approach in~\cite[Section 10.2]{Mor24}, i.e., by implementing a third incidence matrix \(A_C\) for the compressor nodes, in addition to the incidence matrices for the internal nodes \(A_I\) and boundary nodes \(A_B\).
Accordingly, all the port-Hamiltonian pipes and compressors in a given network are stacked in block diagonal form into \(\E(\z), \J(\z), \Rmat(\z)\) and \(\G(\z)\).
Assuming that the nodes are ordered in such a way that the first set of nodes are the boundary nodes, the second set of nodes --- the compressor nodes, and the third set of nodes --- the internal nodes.
Thus, the dynamics of the system are then given by
\begin{equation}
  \begin{split}
    \begin{bmatrix}
      \E(\z) & 0 & 0 & 0 & 0\\
           0 & 0 & 0 & 0 & 0 \\
           0 & 0 & 0 & 0 & 0 \\
           0 & 0 & 0 & 0 & 0 \\
           0 & 0 & 0 & 0 & 0
    \end{bmatrix}
    \begin{bmatrix}
      \dot \z\\
      \dot \mu\\
      \dot \lambda_B\\
      \dot \lambda_C\\
      \dot \lambda_I
    \end{bmatrix}
     &=
     \begin{bmatrix}
      \J(\z) - \Rmat(\z) & \G(\z) & 0 & 0 & 0\\
      -\G(\z)^T & 0 & A_B^T & A_C^T & A_I^T\\
      0 & -A_B & 0 & 0 & 0\\
      0 & -A_C & 0 & 0 & 0\\
      0 & -A_I & 0 & 0 & 0\\
     \end{bmatrix}
     \begin{bmatrix}
      \e(\z)\\
      \mu\\
      \lambda_B\\
      \lambda_C\\
      \lambda_I
    \end{bmatrix}
    +
    \begin{bmatrix}
      0 & 0\\
      0 & 0\\
      I_b & 0\\
      0 & I_c\\
      0 & 0
    \end{bmatrix}
    \begin{bmatrix}
      \u_B\\
      \u_C
    \end{bmatrix}\\
    \y &=
    \begin{bmatrix}
      0 & 0 & I_b & 0 & 0\\
      0 & 0 & 0 & I_c & 0
    \end{bmatrix}
    \begin{bmatrix}
      \e(\z)\\
      \mu\\
      \lambda_B\\
      \lambda_C\\
      \lambda_I
    \end{bmatrix}
    =
    \begin{bmatrix}
      \lambda_B\\
      \lambda_C\\
    \end{bmatrix}.
  \end{split}
\end{equation}
\paragraph*{Example}
In case of the simple network depicted in Fig.~\ref{fig:compressor_network}, the incidence matrices are
\begin{align*}
  A_B &= \begin{bmatrix}
    1 & 0 & 0 & 0 & 0 & 0 & 0 & 0\\
    0 & 0 & 0 & 0 & 0 & 1 & 0 & 0\\
    0 & 0 & 0 & 0 & 0 & 0 & 0 & 1
  \end{bmatrix},\\
  A_C &= \begin{bmatrix}
    0 & 1 & 0 & 0 & 0 & 0 & 0 & 0\\
    0 & 0 & 1 & 0 & 0 & 0 & 0 & 0
  \end{bmatrix} \text{ and }\\
  A_I &= \begin{bmatrix}
    0 & 0 & 0 & 1 & 1 & 0 & 1 & 0
  \end{bmatrix}.
\end{align*}
The boundary inputs are given as \(\u_B = \begin{bmatrix}
  p_\inlet^{\v_1}\\
  -m_\outlet^{\v_2}\\
  -m_\outlet^{\v_3}\\
\end{bmatrix}\), where \(p_\inlet^{\v_1}\) is the inlet pressure at \(\v_1\) and \(m_\outlet^{\v_i}\) are the outlet momenta at \(v_i\), \(i=1,2\).
The inputs for the compressor, \(\u_C\), are, according to assumption \textbf{AV} and framework \textbf{FC},
\begin{align*}
  \u_C^{\textbf{AV},\textbf{FC}} = \begin{bmatrix}
  \tilde{c}^{-\frac{1}{\kappa}}\\
  \tilde{c}
\end{bmatrix}, \quad
  \u_C^{\textbf{AM},\textbf{FC}} = \begin{bmatrix}
  1\\
  \tilde{c}
\end{bmatrix}, \quad
  \u_C^{\textbf{AV},\textbf{FP}} = \begin{bmatrix}
  p_\outlet^{-\frac{1}{\kappa}}\\
  p_\outlet
\end{bmatrix} \quad \text{and}\quad
  \u_C^{\textbf{AM},\textbf{FP}} = \begin{bmatrix}
  1\\
  p_\outlet
\end{bmatrix}.
\end{align*}
\section{Numerical Experiments}

As an indicative experiment, the widely benchmarked Yamal-Europe pipeline testcase is considered for its simplicity, similar to various articles in literature~\cite{HimGB21,Cha09}.
The gas flow over a day is computed using the port-Hamiltonian model across the \qty{363}{\kilo\metre}-long, \qty{1.422}{\metre}-diameter pipeline with a minor modification.
A compressor with a variable configuration is placed exactly midway along the pipeline (Fig~\ref{fig:compressor}).
The gas is supplied into the left pipe consistently with a temperature of \qty{276.25}{\kelvin} and pressure of \qty{80}{\bar}.
A variable flow rate is chosen for gas extraction at the end of the right pipe at every \qty{6}{\hour} period during the day.
The pipe friction factor is considered a constant value of \num[exponent-mode = scientific]{0.0018}.
The fluid is considered ideal with a gas constant of \qty{530}{\joule\per\kilo\gram\per\mole} and an isentropic exponent of \num{1.4}.

The numerical implementation of the model utilizes the finite difference method in space.
The pipe is discretized in space into 32 equal intervals and a steady state simulation is performed for the initial boundary conditions ignoring the time-derivative terms in the port-Hamiltonian model to serve as initial conditions for the transient case.
The transient simulation is subsequently computed for the dynamic discontinuous boundary conditions using an implicit-midpoint integration with a timestep of \qty{100}{\second}.
The numerical scheme is particularly chosen for its numerical stability~\cite{Mor24}.

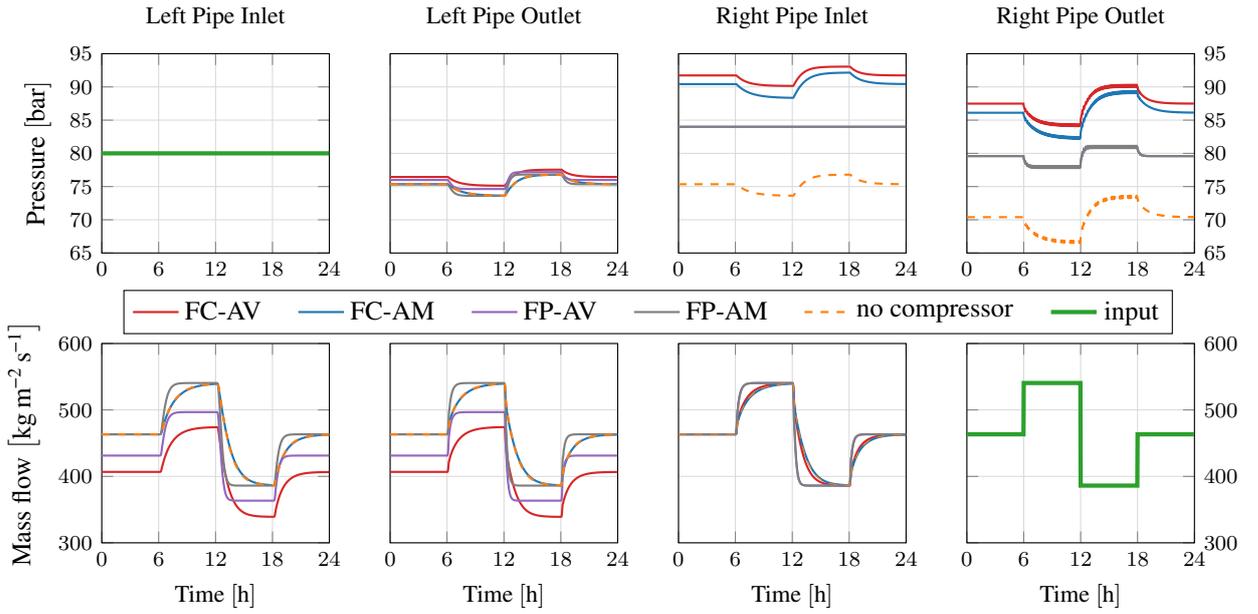
\begin{figure}
  \centering
	\tikzexternalenable%
	\tikzsetnextfilename{results_plot_compressor}%
	\filemodCmp{tikz/results_plot_compressor.tex}{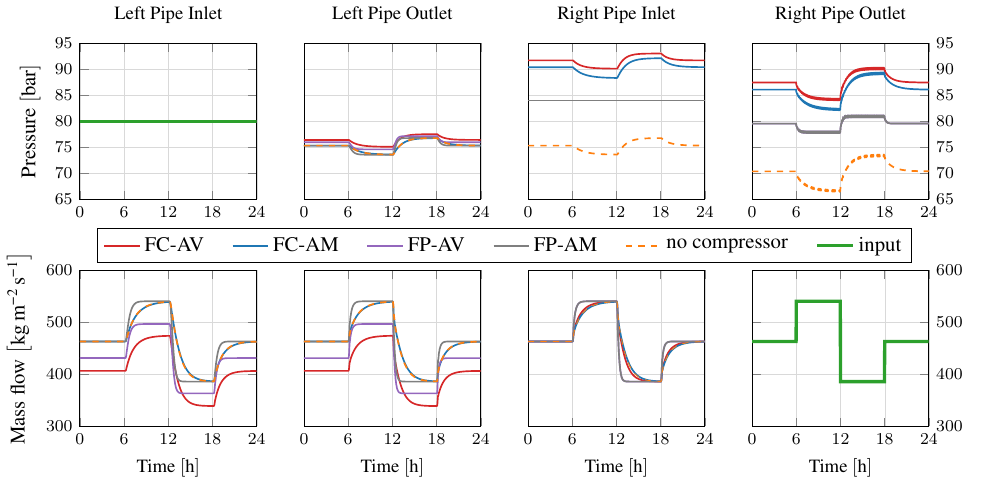}%
	{\tikzset{external/remake next}}{}%
	\pgfplotstableread[col sep=space]{tikz/results_plot_compressor/results_pH_fcav_pipe0.csv}{\fcavL}
\pgfplotstableread[col sep=space]{tikz/results_plot_compressor/results_pH_fcav_pipe1.csv}{\fcavR}

\pgfplotstableread[col sep=space]{tikz/results_plot_compressor/results_pH_fcam_pipe0.csv}{\fcamL}
\pgfplotstableread[col sep=space]{tikz/results_plot_compressor/results_pH_fcam_pipe1.csv}{\fcamR}

\pgfplotstableread[col sep=space]{tikz/results_plot_compressor/results_pH_fpav_pipe0.csv}{\fpavL}
\pgfplotstableread[col sep=space]{tikz/results_plot_compressor/results_pH_fpav_pipe1.csv}{\fpavR}

\pgfplotstableread[col sep=space]{tikz/results_plot_compressor/results_pH_fpam_pipe0.csv}{\fpamL}
\pgfplotstableread[col sep=space]{tikz/results_plot_compressor/results_pH_fpam_pipe1.csv}{\fpamR}

\pgfplotstableread[col sep=space]{tikz/results_plot_compressor/results_pH_nocompressor_pipe0.csv}{\fnocL}
\pgfplotstableread[col sep=space]{tikz/results_plot_compressor/results_pH_nocompressor_pipe1.csv}{\fnocR}

\pgfplotsset{
  every axis/.append style={
    title style={font=\small},
    label style={font=\small},
    tick label style={font=\footnotesize}
  },
  unit code/.code 2 args={\si{#1#2}}
}

\begin{tikzpicture}
  \begin{groupplot}[
    group style={
        group size=4 by 2,
        vertical sep=1.2 cm,
        horizontal sep=0.8cm
    },
    height=0.24\textwidth,
    width=0.26\textwidth,
    xmin=0, xmax=24,
    xtick distance = 6,
    grid=both,
    grid style={line width=.1pt, draw=gray!30},
    no markers, ylabel near ticks,
  ]
    \nextgroupplot[
        title  = Left Pipe Inlet,
        ylabel = {Pressure},
        y unit = {\bar},
        ymin=65, ymax=95, 
        ytick={65, 70, 75, 80, 85, 90, 95},
    ]
      \addplot+[tab20_green, ultra thick] table[x=time, y=inletPressure] {\fcavL};
    \nextgroupplot[
        title  = Left Pipe Outlet,
        ymajorticks = false,
        ymin=65, ymax=95, 
        ytick={65, 70, 75, 80, 85, 90, 95},
        cycle list = {tab20_red, tab20_blue, tab20_purple, tab20_gray, tab20_orange},
        legend to name={CommonLegend},
        legend style={legend columns=6},
        legend cell align={left},
        legend style={/tikz/every even column/.append style={column sep=0.4cm}},
    ]
      \addplot+[thick] table[x=time, y=outletPressure] {\fcavL}; \addlegendentry{FC-AV}
      \addplot+[thick] table[x=time, y=outletPressure] {\fcamL}; \addlegendentry{FC-AM}
      \addplot+[thick] table[x=time, y=outletPressure] {\fpavL}; \addlegendentry{FP-AV}
      \addplot+[thick] table[x=time, y=outletPressure] {\fpamL}; \addlegendentry{FP-AM}
      \addplot+[thick, dashed] table[x=time, y=outletPressure] {\fnocL}; \addlegendentry{no compressor}
      \addlegendimage{tab20_green, ultra thick} \addlegendentry{input}
    \nextgroupplot[
        title  = Right Pipe Inlet,
        ymajorticks = false,
        ymin=65, ymax=95, 
        ytick={65, 70, 75, 80, 85, 90, 95},
        cycle list = {tab20_red, tab20_blue, tab20_purple, tab20_gray, tab20_orange},
    ]
      \addplot+[thick] table[x=time, y=inletPressure] {\fcavR};
      \addplot+[thick] table[x=time, y=inletPressure] {\fcamR};
      \addplot+[thick] table[x=time, y=inletPressure] {\fpavR};
      \addplot+[thick] table[x=time, y=inletPressure] {\fpamR};
      \addplot+[thick, dashed] table[x=time, y=inletPressure] {\fnocR};
    \nextgroupplot[
      title = Right Pipe Outlet,
      yticklabel pos=right,
      ymin=65, ymax=95, 
      ytick={65, 70, 75, 80, 85, 90, 95},
      cycle list = {tab20_red, tab20_blue, tab20_purple, tab20_gray, tab20_orange},
    ]
      \addplot+[thick] table[x=time, y=outletPressure] {\fcavR};
      \addplot+[thick] table[x=time, y=outletPressure] {\fcamR};
      \addplot+[thick] table[x=time, y=outletPressure] {\fpavR};
      \addplot+[thick] table[x=time, y=outletPressure] {\fpamR};
      \addplot+[thick, dashed] table[x=time, y=outletPressure] {\fnocR};
    \nextgroupplot[
        cycle list = {tab20_red, blue},
        xlabel = {Time}, x unit = {\hour},
        ylabel = {Mass flow}, y unit = {\kilogram\per\meter\squared\per\second},
        ymin=300, ymax=600,
        cycle list = {tab20_red, tab20_blue, tab20_purple, tab20_gray, tab20_orange},
    ]
      \addplot+[thick] table[x=time, y=inletMomentum] {\fcavL};
      \addplot+[thick] table[x=time, y=inletMomentum] {\fcamL};
      \addplot+[thick] table[x=time, y=inletMomentum] {\fpavL};
      \addplot+[thick] table[x=time, y=inletMomentum] {\fpamL};
      \addplot+[thick, dashed] table[x=time, y=inletMomentum] {\fnocL};
    \nextgroupplot[
        xlabel = {Time},
        x unit = {\hour},
        ymajorticks = false,
        ymin=300, ymax=600,
        cycle list = {tab20_red, tab20_blue, tab20_purple, tab20_gray, tab20_orange},
    ]
      \addplot+[thick] table[x=time, y=outletMomentum] {\fcavL};
      \addplot+[thick] table[x=time, y=outletMomentum] {\fcamL};
      \addplot+[thick] table[x=time, y=outletMomentum] {\fpavL};
      \addplot+[thick] table[x=time, y=outletMomentum] {\fpamL};
      \addplot+[thick, dashed] table[x=time, y=outletMomentum] {\fnocL};
    \nextgroupplot[
      xlabel = {Time},
      x unit = {\hour},
      ymajorticks = false,
      ymin=300, ymax=600,
      cycle list = {tab20_red, tab20_blue, tab20_purple, tab20_gray, tab20_orange},
    ]
      \addplot+[thick] table[x=time, y=inletMomentum] {\fcavR};
      \addplot+[thick] table[x=time, y=inletMomentum] {\fcamR};
      \addplot+[thick] table[x=time, y=inletMomentum] {\fpavR};
      \addplot+[thick] table[x=time, y=inletMomentum] {\fpamR};
      \addplot+[dashed] table[x=time, y=inletMomentum] {\fnocR};
  \nextgroupplot[
        xlabel = {Time},
        x unit = {\hour},
        yticklabel pos=right,
        ymin=300, ymax=600,
    ]
      \addplot+[tab20_green, ultra thick] table[x=time, y=outletMomentum] {\fcavR};
  \end{groupplot}
  \path (group c2r2.north east) -- node[above]{\ref{CommonLegend}} (group c3r2.north west);
\end{tikzpicture}%
	\tikzexternaldisable%

  \caption{%
    Pressure and Momentum values at each of the nodes for four different models of compressor placed between two pipes.
  }
  \label{fig:compresor_comparison}
\end{figure}

The simulation is performed for compressor models where the specified parameter is either compression ratio (FC) or output pressure (FP).
For the FC framework, the compression ratio is \num{1.2}, whereas for the FP framework, a constant compressor pressure of \qty{84}{\bar} is specified.
The two different frameworks are paired with either a constant velocity assumption (AV)  or a constant momentum assumption (AM).
Fig~\ref{fig:compresor_comparison} shows the result of pressure and momentum values at the inlet and outlet of the two pipes across the compressor.
The results are constrasted against values obtained when there is no compressor present between the pipes.
It is clearly observed that for FC compressors, the inlet pressure at the right pipe is scaled by the compression ratio factor relative to the outlet pressure at the left pipe.
Additionally for FP compressors, the gas into the compressor maintains a fixed outlet pressure.
The AM assumption ensures a constant momentum across the compressor while AV assumption allows variations.
The choice of the right compressor model and assumption is seen to have significant effect on the pressure and momentum fields within the adjacent pipes.

\section{Conclusion}

The objective of this investigation was to explore the inclusion of compressor models in the port-Hamiltonian framework for gas networks.
Four different models of compressors were considered with a pairwise combination of two specifications, either a compression ratio or output pressure; and two assumptions, constant velocity or constant momentum.
The compressor models were imbued into an isothermal port-Hamiltonian model for gas networks by linking the nodes across it through ports.
The numerical implementation of the model was performed and the results displayed the validity of the proposed model.
The study serves as a preliminary outlook into developing an integrated energy-based model for gas networks.

\section*{Software Material}%
  The scripts used in this article are developed using C++ programming language and containerized for reproducibility.
  Eigen library~\cite{eigen} is used to handle operators and linear algebra, and Ceres solver~\cite{ceres} is used to solve the system of non-linear equations.
  The record of the code along with execution instructions are available with {doi:10.5281/zenodo.11387852}~\cite{NayB24}.
\section*{Acknowledgments}%
  The authors thank Riccardo Morandin and Karim Cherifi for helpful discussions.
\vspace{\baselineskip}
\bibliographystyle{plainurl}
\bibliography{literature}
\end{document}